\documentclass[a4paper]{amsart}
\usepackage{amssymb} 
\usepackage[T1]{fontenc}
\usepackage[latin1]{inputenc}
\usepackage{amsfonts}
\usepackage{amsxtra}
\usepackage{ae}
\usepackage{pdfsync}
\usepackage[all]{xy}
\usepackage{enumerate}

\include{diagram}

\newcommand{\Poly}{\mathcal{O}}
\newcommand*{\ket}{\rangle}
\newcommand*{\bra}{\langle}
\newcommand*{\ad}{\mathsf{ad}}

\newcommand*{\E}{\mathcal{E}}

\renewcommand*{\H}{\mathcal{H}}

\newcommand*{\T}{\mathcal{T}}

\newcommand*{\CI}{\mathcal{CI}}

\renewcommand*{\max}{\mathsf{f}}
\newcommand*{\red}{\mathsf{r}}
\newcommand*{\Alg}{$-$\mathsf{Alg}}

\newcommand*{\alg}{\mathsf{alg}}

\newcommand*{\KH}{\mathbb{K}}
\newcommand*{\LH}{\mathbb{L}}

\newcommand*{\Irr}{\mathsf{Irr}}

\newcommand*{\cop}{\mathsf{cop}}

\DeclareMathOperator{\Qut}{Qut}

\DeclareMathOperator{\Aut}{Aut}

\DeclareMathOperator{\coker}{coker}
\DeclareMathOperator{\ind}{ind}
\DeclareMathOperator{\res}{res}

\DeclareMathOperator{\im}{im}

\DeclareMathOperator{\id}{id}

\numberwithin{equation}{section}
\theoremstyle{change}
\newtheorem{theorem}{Theorem}[section]
\newtheorem{prop}[theorem]{Proposition}
\newtheorem{lemma}[theorem]{Lemma}
\newtheorem{cor}[theorem]{Corollary}
\newtheorem{definition}[theorem]{Definition}

\begin{document}

\title[Quantum automorphism groups]{On the structure of quantum automorphism groups}

\author{Christian Voigt}
\address{School of Mathematics and Statistics \\
         University of Glasgow \\
         15 University Gardens \\
         Glasgow G12 8QW \\
         United Kingdom 
}
\email{christian.voigt@glasgow.ac.uk}

\thanks{This work was supported by the Engineering and Physical Sciences Research Council Grant EP/L013916/1.}

\subjclass[2000]{19D55, 81R50}

\maketitle

\begin{abstract}
We compute the $ K $-theory of quantum automorphism groups of finite dimensional $ C^* $-algebras in the sense of Wang. The results show in 
particular that the $ C^* $-algebras of functions on the quantum permutation groups $ S_n^+ $ are pairwise non-isomorphic 
for different values of $ n $. \\
Along the way we discuss some general facts regarding torsion in discrete quantum groups. In fact, the duals of quantum automorphism 
groups are the most basic examples of discrete quantum groups exhibiting genuine quantum torsion phenomena. 
\end{abstract} 

\section{Introduction}

Quantum automorphism groups were introduced by Wang \cite{Wangqsymmetry} in his study of noncommutative symmetries of finite dimensional $ C^* $-algebras. 
These quantum groups are quite different from $ q $-deformations of compact Lie groups, and interestingly, they appear naturally 
in a variety of contexts, including combinatorics and free probability, see for instance \cite{BBCsurvey}, \cite{BCSdefinetti}. 
The $ C^* $-algebraic properties of quantum automorphism groups were studied by Brannan \cite{Brannanquantumautomorphism}, revealing 
various similarities with free group $ C^* $-algebras. \\ 
An interesting subclass of quantum automorphism groups is provided by quantum permutation groups. 
Following \cite{BSliberation}, we will write $ S_n^+ $ for the quantum permutation group on $ n $ letters. According to the definition of Wang, the quantum 
group $ S_n^+ $ is the universal compact quantum group acting on the abelian $ C^* $-algebra $ \mathbb{C}^n $. If one replaces $ \mathbb{C}^n $ by 
a general finite dimensional $ C^* $-algebra, one has to add the data of a state and restrict to state-preserving actions in the definition of 
quantum automorphism groups. Indeed, the choice of state is important in various respects. This is illustrated, for instance, by the work of De Rijdt 
and Vander Vennet on monoidal equivalences among quantum automorphism groups \cite{dRV}. \\ 
The aim of the present paper is to compute the $ K $-theory of quantum automorphism groups. Our general strategy follows the ideas in \cite{Voigtbcfo}, 
which in turn build on methods from the Baum-Connes conjecture, formulated in the language of category theory following 
Meyer and Nest \cite{MNtriangulated}. In fact, the main result of \cite{Voigtbcfo} implies rather easily that the 
appropriately defined assembly map for duals of quantum automorphism groups is an isomorphism. The main additional ingredient, discussed further below, 
is the construction of suitable resolutions, entering the left hand side of the assembly map in the framework of \cite{MNtriangulated}. \\ 
The reason why this is more tricky than in \cite{Voigtbcfo} is that quantum automorphism groups have torsion. At first sight, the  
presence of torsion may appear surprising because these quantum groups behave like free groups in 
many respects. Indeed, the way in which torsion enters the picture is different from what happens for classical discrete groups. 
Therefore quantum automorphism groups provide an interesting class of examples also from a conceptual point of view. Indeed, a better 
understanding of quantum torsion seems crucial in order to go beyond the class of quantum groups studied 
in the spirit of the Baum-Connes conjecture so far \cite{MNcompact}, \cite{Voigtbcfo}, \cite{VVfreeu}. 
We have therefore included some basic considerations on torsion in discrete quantum groups in this paper. \\ 
From our computations discussed below one can actually see rather directly the effect of torsion on the level of $ K $-theory. In particular, 
the $ K $-groups of monoidally equivalent quantum automorphism groups can differ quite significantly due to minor
differences in their torsion structure. Our results also have some direct operator algebraic consequences, most notably, they 
imply that the reduced $ C^* $-algebras of functions on quantum permutation groups can be distinguished by $ K $-theory. \\ 
Let us now explain how the paper is organised. In section \ref{secqg} we collect some definitions and facts from 
the theory of compact quantum groups and fix our notation. Section \ref{secqut} contains more specific 
preliminaries on quantum automorphism groups and their actions. In section \ref{sectorsion} 
we collect some basic definitions and facts regarding torsion in discrete quantum groups. In the quantum case, this
is studied most efficiently in terms of ergodic actions of the dual compact quantum groups, and 
our setup generalises naturally previous considerations by Meyer and Nest \cite{MNcompact}, \cite{Meyerhomalg2}. 
Finally, section \ref{seckqaut} contains our main results.  \\
Let us conclude with some remarks on notation. We write $ \LH(\E) $ for the algebra of adjointable operators on a Hilbert module $ \E $. 
Moreover $ \KH(\E) $ denotes the algebra of compact operators. The closed linear span of a subset $ X $ of a Banach space is denoted 
by $ [X] $. Depending on the context, the symbol $ \otimes $ denotes either the tensor product of Hilbert spaces, 
the spatial tensor product of $ C^* $-algebras, or the exterior tensor product of Hilbert modules.

\section{Compact quantum groups} \label{secqg}

In this preliminary section we collect some definitions from the theory of compact quantum groups and 
fix our notation. We will mainly follow the conventions in \cite{Voigtbcfo} as far as general quantum group theory is concerned. \\ 
Let us start with the following definition. 
\begin{definition} \label{defcqg}
A compact quantum group $ G $ is given by a unital Hopf $ C^* $-algebra $ C(G) $, that is, a unital $ C^* $-algebra $ C(G) $ 
together with a unital $ * $-homomorphism $ \Delta: C(G) \rightarrow C(G) \otimes C(G) $, called comultiplication, such that 
$$
(\Delta \otimes \id) \Delta = (\id \otimes \Delta) \Delta 
$$
and 
$$ 
[(C(G) \otimes 1) \Delta(C(G))] = C(G) \otimes C(G) = [(1 \otimes C(G)) \Delta(C(G))]. 
$$
\end{definition} 
For every compact quantum group there exists a Haar state, namely a state $ \phi: C(G) \rightarrow \mathbb{C} $ satisfying the invariance 
conditions $ (\id \otimes \phi)\Delta(f) = \phi(f) 1 = (\phi \otimes \id)\Delta(f) $ for all $ f \in C(G) $. 
The image of $ C(G) $ in the GNS-representation of $ \phi $ is denoted $ C^\red(G) $, and called the reduced 
$ C^* $-algebra of functions on $ G $. We will write $ L^2(G) $ for the GNS-Hilbert space of $ \phi $, and notice that the GNS-representation 
of $ C^\red(G) $ on $ L^2(G) $ is faithful. \\ 
A unitary representation of $ G $ on a Hilbert space $ \H $ is a unitary element $ U \in M(C^\red(G) \otimes \KH(\H)) = \LH(C^\red(G) \otimes \H) $ 
such that $ (\Delta \otimes \id)(U) = U_{13} U_{23} $. In analogy with the classical theory for compact groups, 
every unitary representation of a compact quantum group $ G $ is completely reducible, and irreducible 
representations are finite dimensional. We write $ \Irr(G) $ for the set of equivalence classes of irreducible unitary 
representations of $ G $. The linear span of matrix coefficients of all unitary representations of $ G $ forms a dense Hopf $ * $-algebra 
$ \Poly(G) $ of $ C^\red(G) $. \\ 
The full $ C^* $-algebra $ C^\max(G) $ of functions on $ G $ is the universal $ C^* $-completion of $ \Poly(G) $. It admits a comultiplication 
as well, satisfying the density conditions in definition \ref{defcqg}. 
The quantum group $ G $ can be equivalently described in terms of $ C^\max(G) $ or $ C^\red(G) $, or in fact, 
using $ \Poly(G) $. One says that $ G $ is coamenable if the canonical quotient map $ C^\max(G) \rightarrow C^\red(G) $ 
is an isomorphism. In this case we will simply write again $ C(G) $ for this $ C^* $-algebra. By slight abuse of notation, 
we will also write $ C(G) $ if a statement holds for both $ C^\max(G) $ and $ C^\red(G) $. \\ 
The regular representation of $ G $ is the representation of $ G $ on $ L^2(G) $ corresponding to the 
multiplicative unitary $ W \in M(C^\red(G) \otimes \KH(L^2(G))) $ determined by 
$$
W^*(\Lambda(f) \otimes \Lambda(g)) = (\Lambda \otimes \Lambda)(\Delta(g)(f \otimes 1)), 
$$
where $ \Lambda(f) \in L^2(G) $ is the image of $ f \in C^\max(G) $ under the GNS-map. 
The comultiplication of $ C^\red(G) $ can be recovered from $ W $ by the formula 
$$
\Delta(f) = W^*(1 \otimes f) W. 
$$
One defines the algebra of functions $ C_0(\hat{G}) $ on the dual discrete quantum group $ \hat{G} $ by 
$$
C_0(\hat{G}) = [(\LH(L^2(G))_* \otimes \id)(W)],
$$
together with the comultiplication 
$$
\hat{\Delta}(x) = \hat{W}^* (1 \otimes x) \hat{W}
$$
for $ x \in C_0(\hat{G}) $, where $ \hat{W} = \Sigma W^* \Sigma $. We remark that 
there is no need to distinguish between $ C_0^\max(\hat{G}) $ and $ C_0^\red(\hat{G}) $ in the discrete case. \\ 
Since we are following the conventions of Kustermans and Vaes \cite{KVLCQG}, there is a flip map built into the above definition of $ \hat{\Delta} $, 
so that the comultiplication of $ C_0(\hat{G}) $ corresponds to the opposite multiplication of $ C^\red(G) $. This is a natural choice 
in various contexts, but it is slightly inconvenient when it comes to Takesaki-Takai duality. We will 
write $ \check{G} $ for $ C_0(\hat{G})^\cop $, that is, for the Hopf $ C^* $-algebra $ C_0(\hat{G}) $ equipped with the 
opposite comultiplication $ \hat{\Delta}^\cop = \sigma \hat{\Delta} $, where $ \sigma $ denotes the flip map. By slight abuse 
of terminology, we shall refer to both $ \check{G} $ and $ \hat{G} $ as the dual quantum group of $ G $, but in the sequel 
we will always work with $ \check{G} $ instead of $ \hat{G} $. According to Pontrjagin duality, the double dual of $ G $ in either of the two conventions 
is canonically isomorphic to $ G $. \\ 
An action of a compact quantum group $ G $ on a $ C^* $-algebra $ A $ is a coaction of $ C(G) $ on $ A $, that is, an injective 
nondegenerate $ * $-homomorphism $ \alpha: A \rightarrow M(C(G) \otimes A) $ 
such that $ (\Delta \otimes \id) \alpha = (\id \otimes \alpha) \alpha $ and $ [(C(G) \otimes 1) \alpha(A)] = C(G) \otimes A $. 
In a similar way one defines actions of discrete quantum groups, or in fact arbitrary locally compact quantum groups. We will 
call a $ C^* $-algebra equipped with a coaction of $ C^\red(G) $ a $ G $-$ C^* $-algebra. Moreover we write $ G \Alg $ for the 
category of all $ G $-$ C^* $-algebras and equivariant $ * $-homomorphisms. \\ 
The reduced crossed product $ G \ltimes_\red A $ of a $ G $-$ C^* $-algebra $ A $ is the $ C^* $-algebra 
$$
G \ltimes_\red A = [(C_0(\check{G}) \otimes 1)\alpha(A)]
$$
The crossed product is equipped with a canonical dual action of $ \check{G} $, 
which turns it into a $ \check{G} $-$ C^* $-algebra. Moreover, one has the following analogue of the Takesaki-Takai duality theorem \cite{BSUM}.
\begin{theorem} \label{TTduality}
Let $ G $ be a regular locally compact quantum group and let $ A $ be a $ G $-$ C^* $-algebra. Then there is a natural isomorphism
$$
\check{G} \ltimes_\red G \ltimes_\red A \cong \KH(L^2(G))) \otimes A
$$
of $ G $-$ C^* $-algebras.
\end{theorem}
We will use Takesaki-Takai duality only for discrete and compact quantum groups, and in this setting regularity is automatic. 
At some points we will also use the full crossed product $ G \ltimes_\max A $ of a $ G $-$ C^* $-algebra $ A $, and we refer to \cite{NVpoincare} 
for a review of its definition in terms of its universal property for covariant representations.

\section{Quantum automorphism groups} \label{secqut}

In this section we review some basic definitions and results on quantum automorphism groups of finite dimensional 
$ C^* $-algebras and fix our notation. We refer to \cite{Wangqsymmetry}, \cite{Banicageneric}, \cite{Banicafusscatalan} for more background on quantum automorphism groups. \\ 
Let us start with the definition of the quantum automorphism group of a finite dimensional $ C^* $-algebra $ A $, 
compare \cite{Wangqsymmetry} \cite{Banicageneric}. If $ \mu: A \otimes A \rightarrow A $ denotes the multiplication map then a faithful 
state $ \omega $ on $ A $ is called a $ \delta $-form for $ \delta > 0 $ if $ \mu \mu^* = \delta^2 \id $ with respect to the Hilbert space 
structures on $ A $ and $ A \otimes A $ implemented by the GNS-constructions for $ \omega $ and $ \omega \otimes \omega $, respectively. 
\begin{definition} \label{defqut}
Let $ A $ be a finite dimensional $ C^* $-algebra and let $ \omega $ be a $ \delta $-form on $ A $ for some $ \delta > 0 $. The quantum automorphism 
group $ \Qut(A) = \Qut(A, \omega) $ is the universal compact quantum group acting on $ A $ such that $ \omega $ is preserved. \\
That is, if $ G $ is any compact quantum group together with a coaction $ \delta: A \rightarrow C^\max(G) \otimes A $ 
then there exists a unique morphism of quantum groups $ \iota: G \rightarrow \Qut(A) $ such that the diagram 
$$
\xymatrix{
A \ar@{->}[r]^{\!\!\!\!\!\!\!\!\!\!\!\!\!\!\!\!\!\!\!\!\!\! \alpha} \ar@{->}[rd]_{\delta} & C^\max(\Qut(A)) \otimes A \ar@{->}[d]^{\iota^* \otimes \id} \\ 
& C^\max(G) \otimes A 
}
$$
is commutative and 
$$
(\id \otimes \omega)\alpha(a) = \omega(a) 1 
$$
for all $ a \in A $. Here $ \iota^*: C^\max(\Qut(A)) \rightarrow C^\max(G) $ denotes the homomorphism of Hopf $ C^* $-algebras corresponding 
to the morphism $ \iota $. 
\end{definition} 
A basic result of the theory is that the compact quantum group $ \Qut(A, \omega) $ indeed exists \cite{Wangqsymmetry}, see also 
\cite{Mrozinskiso3deformations}. As indicated in definition \ref{defqut}, we will typically omit the state $ \omega $ from our notation and write 
$ \Qut(A) $ instead of $ \Qut(A,\omega) $, although $ \omega $ is an important part of the data. 
The notation for quantum automorphism groups used in \cite{Wangqsymmetry} is $ A_{aut}(B) = C^\max(\Qut(B)) $. \\ 
We remark that the matrix coefficients of an action of a compact quantum 
group $ G $ on a finite dimensional $ C^* $-algebra $ A $ are contained in the Hopf $ * $-algebra $ \Poly(G) $. 
In particular, any coaction $ \alpha: A \rightarrow C^\max(G) \otimes A $, where $ G $ is a compact quantum group, 
comes from a Hopf algebraic coaction $ \alpha: A \rightarrow \Poly(G) \otimes A $. Using this fact, one can study 
quantum automorphism groups from a more algebraic perspective, see for instance \cite{Mrozinskiso3deformations}. The 
universal $ C^* $-algebras of quantum automorphism groups can be defined explicitly in terms of generators and relations \cite{Wangqsymmetry}. \\ 
Let us take a closer look at the special case of quantum permutation groups. 
By definition, the quantum permutation group $ S_n^+ $ is the quantum automorphism group of $ A = \mathbb{C}^n $ with the 
trace corresponding to the uniform probability measure on $ n $ points, which is a $ \delta $-form for $ \delta = \sqrt{n} $. 
In order to describe $ C^\max(S_n^+) $ let us recall some terminology. If $ B $ is any unital $ * $-algebra then a matrix $ u = (u_{ij}) \in M_n(B) $ 
is called a magic unitary if all entries $ u_{ij} $ are projections, and on each row and column of $ u $ these projections are 
mutually orthogonal and sum up to $ 1 $. Explicitly, this means 
\begin{equation*}
u_{ij}^* = u_{ij} = u_{ij}^2 
\end{equation*}
for all $ 1 \leq i,j \leq n $ and 
\begin{equation*}
\sum_{i = 1}^n u_{ik} = 1, \qquad \sum_{i = 1}^n u_{ki} = 1 
\end{equation*}
for all $ k $. These relations imply in particular that the matrix $ u $ and its transpose $ u^t $ are both unitary. 
\begin{prop}
The full $ C^* $-algebra of functions on the quantum permutation group $ S_n^+ = \Qut(\mathbb{C}^n) $ is 
the universal unital $ C^* $-algebra $ C^\max(S_n^+) $ with generators $ u_{ij} $ for $ 1 \leq i,j \leq n $ 
such that $ u = (u_{ij}) \in M_n(C^\max(S_n^+)) $ is a magic unitary matrix. 
\end{prop} 
The comultiplication $ \Delta: C^\max(S_n^+) \rightarrow C^\max(S_n^+) \otimes C^\max(S_n^+) $ is defined by the formula
$$
\Delta(u_{ij}) = \sum_{k = 1}^n u_{ik} \otimes u_{kj} 
$$
on generators. 
The defining coaction $ \alpha: \mathbb{C}^n \rightarrow C^\max(S_n^+) \otimes \mathbb{C}^n $ is given by 
$$
\alpha(e_i) = \sum_{j = 1}^n u_{ij} \otimes e_j,  
$$
where $ e_1, \dots, e_n $ are the minimal projections in $ \mathbb{C}^n $. \\ 
From the fact that $ S_n^+ $ is the universal compact quantum group acting on $ A = \mathbb{C}^n $ we obtain a morphism of quantum groups 
$ S_n \rightarrow S_n^+ $, that is, a unital $ * $-homomorphism $ C^\max(S_n^+) \rightarrow C(S_n) $ compatible with comultiplications. 
Here $ S_n $ is the symmetric group acting by permutations on $ A $. 
In fact, every character of $ C^\max(S_n^+) $ is induced from a character of $ C(S_n) $, and $ C(S_n) $ is the abelianisation of 
$ C^\max(S_n^+) $. \\ 
The structure of quantum permutation groups $ S_n^+ $ is well-understood for small values of $ n $, 
for the following fact compare \cite{Banicasmallmetric}. 
\begin{lemma} \label{s123plusstructure}
For $ n = 1,2,3 $ the canonical morphism of quantum groups $ S_n \rightarrow S_n^+ $ is an isomorphism. 
\end{lemma} 
Clearly, this means in particular that $ S_n^+ $ is coamenable for these values of $ n $. \\
For $ n = 4 $ the natural morphism $ S_n \rightarrow S_n^+ $ is no longer an isomorphism. 
In fact, the $ C^* $-algebra $ C^\max(S_4^+) $ is infinite dimensional. 
The following result due to Banica and Bichon \cite{BBfourpoints} describes the structure of $ S_4^+ $. 
\begin{theorem} \label{snplus4structure}
There is an isomorphism of quantum groups $ S_4^+ \cong SO_{-1}(3) $, where the latter is obtained 
from $ SU_{-1}(3) $ by making the fundamental matrix orthogonal. The quantum group $ SO_{-1}(3) $ 
is a $ 2 $-cocycle twist of the classical group $ SO(3) $. 
\end{theorem} 
We note that since the classical group $ SO(3) $ is a coamenable compact quantum group, the same holds for its cocycle 
twist $ S_4^+ $, see \cite{Banicafusion}. 
For $ n \geq 5 $ the quantum groups $ S_n^+ $ are no longer coamenable, that is, the reduced $ C^* $-algebras $ C^\red(S_n^+) $ 
fail to be nuclear. \\
Still, the $ C^* $-algebras $ C^\red(S_n^+) $ are exact for all $ n $. This can be shown using the monoidal equivalences among 
quantum automorphism groups \cite{dRV} to be discussed below, by invoking a general observation of Vaes and Vergnioux \cite{VaesVergnioux}, 
namely that exactness of the reduced $ C^* $-algebras of functions on compact quantum groups is preserved under monoidal equivalences, 
compare \cite{Brannanquantumautomorphism}. \\ 
According to lemma \ref{s123plusstructure}, all quantum automorphisms for $ C^* $-algebras of dimension at most $ 3 $ come from classical 
automorphisms. In the sequel we will therefore restrict attention to finite dimensional $ C^* $-algebras of dimension at least $ 4 $. \\ 
In the case of dimension $ 4 $, the only $ C^* $-algebra to consider apart from $ \mathbb{C}^4 $ is $ M_2(\mathbb{C}) $. 
The quantum automorphism group of $ A = \mathbb{C}^4 $ is the quantum permutation group $ S_4^+ $ 
already mentioned above, and for $ M_2(\mathbb{C}) $, the quantum automorphism groups are determined by the following result 
of Soltan \cite{Soltanquantumso3}. 
\begin{theorem} 
Let $ \omega $ be a faithful state on $ M_2(\mathbb{C}) $. The quantum automorphism group $ \Qut(M_2(\mathbb{C}), \omega) $ is isomorphic 
to $ SO_q(3) $ for a unique $ q \in (0,1] $. Here $ SO_q(3) $ is the quantum $ SO(3) $-group of Podle\`s. 
\end{theorem}
We emphasise that the quantum group $ SO_q(3) $ of Podle\`s \cite{Podlesspheres} is different from the quantum group $ SO_{-1}(3) $ appearing 
in theorem \ref{snplus4structure}. \\ 
Let us now return to general quantum automorphism groups. The main 
tool at our disposal are the monoidal equivalences exhibited by De Rijdt and Vander Vennet as follows \cite{dRV}. 
\begin{theorem} \label{fdmonoidaleq}
Let $ A_j $ be finite dimensional $ C^* $-algebras of dimension at least $ 4 $, equipped with $ \delta_j $-forms $ \omega_j $ 
for $ j = 1,2 $, respectively. Then $ \Qut(A_1) $ is monoidally equivalent to $ \Qut(A_2) $ iff $ \delta_1 = \delta_2 $. 
\end{theorem} 
This result will play a crucial role in our analysis, it implies in particular that any quantum automorphism group 
is monoidally equivalent to $ SO_q(3) $ for some $ q \in (0,1] $. 
Theorem \ref{fdmonoidaleq} shows also that the quantum permutation groups $ S_n^+ = \Qut(\mathbb{C}^n) $ are pairwise 
monoidally inequivalent for $ n \geq 4 $.

\section{Torsion in discrete quantum groups} \label{sectorsion} 

In this section we discuss some definitions and facts related to torsion in discrete quantum groups. This will be useful in our analysis 
of the equivariant Kasparov theory of quantum automorphism groups. The study of torsion in duals of compact groups 
has already been carried out in \cite{MNcompact}, and our definitions are by and large motivated from this. \\
Firstly, we have to explain what we mean by torsion. Let $ \check{G} $ be the dual of a discrete quantum group $ G $, 
and assume that $ U \in C^\red(\check{G}) \otimes \KH(V) $ is a unitary representation of $ \check{G} $ on the finite dimensional Hilbert 
space $ V $. Then we obtain a coaction $ \ad_U: \KH(V) \rightarrow C^\red(\check{G}) \otimes \KH(V) $ by the formula 
$$
\ad_U(T) = U^* (\id \otimes T) U, 
$$
which turns $ \KH(\H) $ into a $ \check{G} $-$ C^* $-algebra. Coactions of this form are precisely 
the coactions on simple matrix algebras which are $ \check{G} $-equivariantly Morita equivalent to the trivial coaction on $ \mathbb{C} $. \\
The idea is, roughly, that torsion in $ G $ is encoded in coactions on finite dimensional $ C^* $-algebras which are not of the above form. 
In order to make this precise, recall that a unital $ \check{G} $-$ C^* $-algebra $ B $ with coaction $ \beta: B \rightarrow C(\check{G}) \otimes B $ 
is called ergodic iff its fixed point subalgebra 
$$
B^\beta = \{b \in B \mid \beta(b) = 1 \otimes b \} 
$$
is equal to $ \mathbb{C} $. \\ 
The following definition is motivated by the considerations regarding torsion-free quantum groups in \cite{Meyerhomalg2}. 
\begin{definition} \label{deftorsion}
Let $ G $ be a discrete quantum group. A finite dimensional ergodic $ \check{G} $-$ C^* $-algebra $ B $ is called a torsion coaction 
of $ G $. Such a coaction is called nontrivial if $ B $ is not $ \check{G} $-equivariantly Morita equivalent to the 
trivial $ \check{G} $-$ C^* $-algebra $ \mathbb{C} $. \\ 
The quantum group $ G $ is called torsion-free if $ G $ does not admit any nontrivial torsion coactions. 
\end{definition}
It is straightforward to check that a quantum group $ G $ is torsion-free iff for every finite dimensional $ \check{G} $-$ C^*$-algebra $ B $ there 
are finite dimensional unitary corepresentations $ U_j \in C^\red(\check{G}) \otimes \KH(V_j) $ such that 
$$ 
B \cong \KH(V_1) \oplus \cdots \oplus \KH(V_l) 
$$ 
as $ \check{G} $-$ C^* $-algebras, where each matrix block $ \KH(V_j) $ is equipped with the adjoint coaction $ \ad_{U_j} $ as explained
above. That is, the notion of torsion-freeness in definition \ref{deftorsion} is compatible with the terminology used in 
\cite{Meyerhomalg2}, \cite{Voigtbcfo}. \\
Assume that $ G $ is a discrete quantum group, and let $ Q $ be a Galois object for $ C^*(H)^\cop = C(\check{H}) $ where 
$ H \subset G $ is a finite quantum subgroup. That is, $ Q $ is a $ C^* $-algebra equipped with an ergodic coaction of $ C(\check{H}) $ 
of full quantum multiplicity, compare for instance \cite{BdRV}, \cite{Decommergaloisalg}. Then $ Q $ is a torsion coaction of $ G $ since $ Q $ 
is clearly finite dimensional, and the ergodic action of $ \check{H} $ on $ Q $ is naturally an ergodic action of $ \check{G} $ as well. \\ 
The following proposition shows that this class of coactions exhausts all torsion coactions in the case of classical discrete groups.  
\begin{prop} \label{discretegrouptorsion}
Let $ G $ be a discrete group. Then every torsion coaction of $ G $ is $ \check{G} $-equivariantly Morita equivalent to 
a $ \check{G} $-$ C^* $-algebra of the form $ C^*_\omega(H) $ for some finite subgroup $ H \subset G $ and a normalised 
cocycle $ \omega \in Z^2(H, U(1)) $. 
\end{prop} 
\proof Let $ \beta: B \rightarrow C^*_\red(G) \otimes B $ be a torsion coaction and consider the corresponding spectral decomposition 
$$
B = \bigoplus_{s \in G} B_s 
$$
of $ B $, where $ B_s = \{b \in B \mid \beta(b) = s \otimes b \} $. 
Ergodicity means that the component of the identity element $ e \in G $ is $ B_e = \mathbb{C} $. 
Observe next that if $ b \in B_s $ then $ b^* \in B_s^{-1} $ and $ b^* b \in B_e $, 
so that $ b^* b = \lambda > 0 $ is invertible. It follows that
all spectral subspaces are one-dimensional and spanned by invertible elements. \\
For every $ s \in G $ such that $ B_s $ is nonzero we may choose a unitary $ \delta_s \in B_s $, and without loss of generality 
we may pick $ \delta_e = 1 $. Clearly, the elements $ s \in G $ such that $ B_s \neq 0 $ form a finite subgroup $ H $ of $ G $. 
Moreover we have $ \delta_s \delta_t = \omega(s,t) \delta_{st} $ for $ \omega(s,t) \in U(1) $. It is straightforward to check that this yields 
a normalised $ 2 $-cocycle $ \omega \in Z^2(H, U(1)) $, and we conclude that $ B $ is isomorphic to $ C^*_\omega(H) $, the twisted group $ C^* $-algebra 
of $ H $. \qed \\
As a corollary of proposition \ref{discretegrouptorsion} we see in particular that the notion of torsion-freeness 
introduced in definition \ref{deftorsion} agrees with the usual terminology for discrete groups. That is, a discrete 
group $ G $ is torsion-free in the sense of definition \ref{deftorsion} iff it has no nontrivial finite subgroups. \\ 
In a sense, the torsion coactions obtained from Galois objects for finite quantum subgroups are the most 
obvious examples of torsion coactions. If one goes beyond classical discrete groups then more exotic torsion coactions can appear. \\ 
In particular, as shown in \cite{MNcompact}, this already happens for duals of compact groups. If $ G $ is a compact 
group then a torsion coaction of the dual discrete quantum group $ \check{G} $ is nothing but an ergodic action of $ G $ 
on a finite dimensional $ C^* $-algebra $ B $. Such actions factorise over a Lie group quotient $ K $ of $ G $ 
because the group $ \Aut(B) $ of $ * $-automorphisms of $ B $ is a compact Lie group. Moreover, ergodicity implies that $ B $ must consist of matrix blocks 
of the same size. That is, we have $ B \cong M_k(\mathbb{C})^{\oplus n} $ for some $ k,n \in \mathbb{N} $. 
The subgroup $ L \subset K $ preserving a fixed matrix block of $ B $ contains the connected component $ K_{(0)} $ of $ K $, and $ B $ 
is isomorphic to the induced $ C^* $-algebra $ \ind_L^K(M_k(\mathbb{C})) $ corresponding to the action of $ L $ on $ M_k(\mathbb{C}) $. 
In other words, we see in particular that all proper homogeneous $ \check{G} $-$ C^* $-algebras considered in \cite{MNcompact}
arise from torsion coactions of $ \check{G} $. \\ 
It is sometimes useful to distinguish between two basic types of torsion coactions. 
Let us say that a projective torsion coaction of a discrete quantum group $ G $ is a torsion coaction whose 
underlying $ C^* $-algebra is simple, and which is not equivariantly Morita equivalent to $ \mathbb{C} $. Moreover, let us refer to all 
torsion coactions on non-simple $ C^* $-algebras as permutation torsion coactions. 
Clearly, the quantum group $ G $ is torsion-free if it has neither projective torsion nor permutation torsion. 
The terminology is motivated from the fact that one sourse of permutation torsion for $ G $ comes from finite quantum subgroups and their 
regular coactions. Similarly, projective torsion arises from projective representations of $ \check{G} $. \\ 
Roughly speaking, the following lemma shows that permutation torsion for $ G $ is related to the connectedness of $ \check{G} $. 
\begin{lemma} 
The dual of a compact group $ G $ has no permutation torsion iff $ G $ is connected. 
\end{lemma} 
\proof If $ G $ is connected, then every action of $ G $ on 
$$
B = M_{k_1}(\mathbb{C}) \oplus \cdots \oplus M_{k_n}(\mathbb{C}) 
$$
must preserve the individual matrix blocks, and ergodicity implies that $ n = 1 $. 
Conversely, assume that $ G $ is not connected. Then the quotient $ G/G_{(0)} $ of $ G $ by its connected component 
is a nontrivial profinite group. If $ F $ is a nontrivial finite quotient of $ G/G_{(0)} $, 
then the permutation action of $ F $ on the commutative $ C^* $-algebra $ C(F) $ induces an ergodic action of $ G $ 
via the quotient maps $ G \rightarrow G/G_{(0)} \rightarrow F $. Hence $ G $ has permutation torsion in this case. \qed \\
Assume that $ G $ and $ H $ are discrete quantum groups such that $ \check{G} $ and $ \check{H} $ are monoidally equivalent. 
Then the general correspondence between actions of $ \check{G} $ and $ \check{H} $ shows that torsion coactions of $ G $ and $ H $ are 
in a bijective correspondence \cite{dRV}, \cite{Voigtbcfo}. 
For duals of quantum automorphism groups it is therefore quite easy to determine all torsion coactions up to equivariant Morita 
equivalence. 
\begin{lemma} \label{quttorsion} 
Let $ \Qut(A, \omega) $ be the quantum automorphism group of a $ C^* $-algebra $ A $ of dimension at least $ 4 $ with 
respect to the $ \delta $-form $ \omega $. Then the trivial coaction on $ \mathbb{C} $ and the defining coaction on $ A $ are the only 
torsion coactions of the dual of $ \Qut(A, \omega) $ up to equivariant Morita equivalence. 
\end{lemma} 
\proof According to the results in \cite{dRV}, the quantum group $ \Qut(A,\omega) $ is monoidally equivalent to $ H = SO_q(3) $ for some $ q \in (0,1] $. 
If we write $ G = SU_q(2) $, then we have $ C(H) \subset C(G) $, and if $ B $ is a finite dimensional ergodic $ H $-$ C^* $-algebra, it is 
also naturally an ergodic $ G $-$ C^* $-algebra. 
According to \cite{Voigtbcfo}, the quantum group $ \check{G} $ is torsion-free. It follows that $ B \cong \KH(V(n)) $ for 
some $ n \in \frac{1}{2} \mathbb{N}_0 \cong \Irr(SU_q(2)) $. If $ n $ is an integer then $ B $ is $ H $-equivariantly Morita 
equivalent to $ \mathbb{C} $. In this case the corresponding ergodic action of $ \Qut(A, \omega) $ 
is Morita equivalent to the trivial action on $ \mathbb{C} $. 
Otherwise $ B $ is $ H $-equivariantly Morita equivalent to $ \KH(V(1/2)) $. 
The $ H $-$ C^* $-algebra $ \KH(V(1/2)) $ corresponds to $ A $ under the monoidal equivalence, because up to isomorphism it is the only 
ergodic $ H $-$ C^* $-algebra with the correct spectral decomposition. \qed \\ 
In particular, lemma \ref{quttorsion} shows that projective torsion may be turned into permutation torsion
under monoidal equivalence, and vice versa. \\ 
For the applications to discrete quantum groups we have to study the crossed products of torsion coactions, 
and it will be convenient to use the following terminology, compare \cite{MNcompact}. 
\begin{definition} 
Let $ G $ be a discrete quantum group. A $ G $-$ C^* $-algebra is called proper almost homogeneous if it is $ G $-equivariantly Morita 
equivalent to the crossed product of some torsion coaction of $ G $. 
\end{definition} 
Notice that we do not need to distinguish between reduced or maximal crossed products here since the dual of $ G $ is compact. \\ 
A guiding example for proper almost homogeneous algebras arises from the torsion coaction $ B = C(\check{H}) $ for some finite quantum 
subgroup $ H \subset G $. In this case the proper almost homogeneous algebra $ G \ltimes B $ is $ G $-equivariantly Morita equivalent to $ C_0(G/H) $, 
which should be viewed as the prototypical example of a proper homogeneous action of $ G $. \\
According to proposition \ref{discretegrouptorsion}, proper almost homogeneous actions of classical discrete groups are indeed 
essentially determined by homogeneous spaces $ G/H $ where $ H $ is finite. Notice that in the presence of a nontrivial cocycle $ \omega \in Z^2(H, U(1)) $, 
the crossed product $ P = \check{G} \ltimes C^*_\omega(H) $ is a proper $ G $-$ C^* $-algebra over $ C_0(G/H) $ in the sense 
of Kasparov. \\ 
Assume that $ G $ be a locally compact quantum group and let $ A $ be a $ G $-$ C^* $-algebra. Let us say that $ G $ 
acts amenably on $ A $ if the canonical quotient map $ G \ltimes_\max A \rightarrow G \ltimes_\red A $ is an isomorphism. 
In particular, according to this terminology, $ G $ is amenable in the sense that $ C^*_\max(G) \cong C^*_\red(G) $ iff $ G $ acts 
amenably on $ \mathbb{C} $. 
\begin{lemma} \label{properamenable}
Let $ G $ be a discrete quantum group. If $ P $ is a proper almost homogeneous $ G $-$ C^* $-algebra then $ G $ acts amenably on $ P $. 
\end{lemma}  
\proof Without loss of generality we may assume that $ P = \check{G} \ltimes B $ for a finite dimensional $ C^* $-algebra $ B $ with an ergodic 
action of $ \check{G} $. It is enough to observe that the algebraic crossed product $ G \ltimes_\alg \check{G} \ltimes_\alg B $, 
taken in the framework of algebraic quantum groups \cite{vDadvances}, is dense inside both $ G \ltimes_\max P $ and $ G \ltimes_\red P $. 
Indeed, the algebraic crossed product is isomorphic to the algebraic 
tensor product of $ B $ with an algebra of possibly infinite matrices, and the $ C^* $-norm on such an 
algebra is uniquely determined. Therefore the quotient map $ G \ltimes_\max P \rightarrow G \ltimes_\red P $ is an isomorphism. \qed

\section{The $ K $-theory of quantum automorphism groups} \label{seckqaut}

In this section we compute the $ K $-theory of quantum automorphism groups. The basic strategy is the same 
as for free quantum groups \cite{Voigtbcfo}, \cite{VVfreeu}, namely, in a first step it will be shown that 
quantum automorphism groups satisfy a strong form of the Baum-Connes conjecture. The second step consists of the actual 
computation, essentially this amounts to computing the left hand side of the assembly map. \\ 
Let $ G = \Qut(A, \omega) $ be the quantum automorphism group of a finite dimensional $ C^* $-algebra $ A $ with respect to 
a $ \delta $-form $ \omega $. In the sequel we will always assume that the dimension of $ A $ is at least $ 4 $. 
We will also fix $ q \in (0,1] $ such that $ H = SO_q(3) $ is monoidally equivalent to $ G $, see \cite{dRV}. \\ 
Firstly, we have to analyse the structure of the equivariant $ KK $-theory of $ G $. For this we need the language of 
triangulated categories, compare \cite{MNtriangulated}, \cite{NVpoincare}, \cite{Voigtbcfo}. More precisely, we consider the 
category $ KK^G $ with objects all separable $ G $-$ C^* $-algebras, and the morphism set between objects $ B $ and $ C $ is given 
by the equivariant Kasparov group $ KK^G(B,C) $. Composition of morphisms is given by the Kasparov product. The category $ KK^G $ is triangulated with 
translation automorphism $ \Sigma: KK^G \rightarrow KK^G $ given by the suspension $ \Sigma B = C_0(\mathbb{R}, B) $ of a $ G $-$ C^* $-algebra $ B $.
The exact triangles are all diagrams isomorphic to diagrams of the form
$$
\xymatrix{
\Sigma C  \;\; \ar@{->}[r] & C_f \ar@{->}[r] & B \ar@{->}[r]^f & C
}
$$
where $ C_f $ denotes the mapping cone of an equivariant $ * $-homomorphism $ f: B \rightarrow C $. 
For more information we refer to \cite{MNtriangulated}, \cite{NVpoincare}, \cite{Voigtbcfo}. \\ 
Let $ \T_G \subset KK^G $ be the full subcategory consisting of all trivial $ G $-$ C^* $-algebras, and all 
$ G $-$ C^* $-algebras of the form $ A \otimes C $, equipped with the the action coming from the defining action of $ G $ on the first 
tensor factor. These actions should be thought of as crossed products of compactly induced actions for the dual discrete quantum group $ \check{G} $. 
We write $ \bra \T_G \ket $ for the localising subcategory of $ KK^G $ generated by $ \T_G $. 
Accordingly, we let $ \bra \CI_{\check{G}} \ket \subset KK^{\check{G}} $ be the full subcategory corresponding to $ \bra \T_G \ket $ 
under Baaj-Skandalis duality \cite{BSUM}, that is, under the equivalence $ KK^G \rightarrow KK^{\check{G}} $ 
of triangulated categories given by taking reduced crossed products. In the sequel we will however avoid working in $ KK^{\check{G}} $ 
for most of the time, because the constructions are somewhat clearer in the compact picture. \\ 
Since $ H = SO_q(3) $ is monoidally equivalent to $ G $, we have an equivalence of triangulated 
categories between $ KK^H $ and $ KK^G $, see \cite{Voigtbcfo}. Notice that the category $ \T_G $ corresponds 
to $ \T_H $ under this equivalence. Therefore, it largely suffices to study the 
structure of $ KK^H $. \\
Recall that the quantum group $ H = SO_q(3) $ is closely related to $ K = SU_q(2) $ 
since $ C(SO_q(3)) \subset C(SU_q(2)) $. We will identify the set $ \Irr(SU_q(2)) $ of equivalence classes of 
irreducible representations of $ SU_q(2) $ with $ \frac{1}{2} \mathbb{N}_0 $, and refer to the irreducible 
representation $ V(n) $ of $ SU_q(2) $ corresponding to $ n \in \frac{1}{2} \mathbb{N}_0 $ as the representation of spin $ n $. 
In the Peter-Weyl picture, the algebra $ C(SO_q(3)) $ is the norm closure of the space of matrix elements of all integral spin 
representations. Accordingly, the set $ \Irr(SO_q(3)) $ of irreducible representations of $ SO_q(3) $ identifies with $ \mathbb{N}_0 $.
We shall write $ C^*_\omega(SO_q(3)) $ for the quotient of $ C^*(SU_q(2)) $ corresponding to the remainig 
part of $ \Irr(SU_q(2)) $. That is, $ C^*_\omega(SO_q(3)) $ is the norm closure of the matrix algebras associated to 
representations of half-integral spin. By construction, we have a direct sum decomposition 
$$
C^*(SU_q(2)) = C^*(SO_q(3)) \oplus C^*_\omega(SO_q(3)), 
$$
and this decomposition is compatible with the canonical coactions of $ C^*(SO_q(3)) $ 
induced from the comultiplication of $ C^*(SU_q(2)) $. 
\begin{lemma} \label{soq3morita}
There is a $ C^*(SO_q(3)) $-colinear Morita equivalence 
$$ 
C^*_\omega(SO_q(3)) \sim SO_q(3) \ltimes \KH(V(1/2)) 
$$ 
where $ V(1/2) $ is the defining representation of $ SU_q(2) $. 
\end{lemma} 
\proof We simply have to cut down the $ C^*(SU_q(2)) $-colinear Morita equivalence 
between $ C^*(SU_q(2)) $ and $ SU_q(2) \ltimes \KH(V(1/2)) $ implemented by the imprimitivity bimodule 
$$ 
C^*(SU_q(2)) \ltimes V(1/2) = [(C^*(SU_q(2)) \otimes 1) \lambda(V(1/2))] \subset \LH(L^2(SU_q(2)) \otimes V(1/2)). 
$$
More precisely, split the Hilbert space $ L^2(SU_q(2)) = \H_0 \oplus \H_1 $ into the direct sum of $ \H_0 = L^2(SO_q(3)) $ and its orthogonal 
complement $ \H_1 $. Then we obtain a Hilbert $ C^*_\omega(SO_q(3)) $-module 
$$ 
[(C^*(SO_q(3)) \otimes 1) \lambda(V(1/2)) (C^*_\omega(SO_q(3)) \otimes 1)] \subset \LH(\H_1, \H_0) \otimes V(1/2), 
$$ 
and it is straightforward to check that this module implements a $ C^*(SO_q(3)) $-colinear Morita equivalence between 
$ C^*_\omega(SO_q(3)) $ and $ SO_q(3) \ltimes \KH(V(1/2)) $. \qed \\ 
The advantage of $ SO_q(3) \ltimes \KH(V(1/2)) $ over $ C^*_\omega(SO_q(3)) $ is that it is easy to check what happens to the 
former under monoidal equivalences. Indeed, according to proposition \ref{quttorsion}, the $ SO_q(3) $-$ C^* $-algebra $ \KH(V(1/2)) $ 
corresponds to the defining action of $ G = \Qut(A) $ on $ A $. \\ 
Recall that a discrete quantum group $ F $ is called $ K $-amenable if the canonical 
map $ F \ltimes_\max B \rightarrow F \ltimes_\red B $ is an isomorphism in $ KK $ for any $ F $-$ C^* $-algebra $ B $, 
compare \cite{Vergniouxkam}. \\ 
The following theorem shows that duals of quantum automorphism groups satisfy the Baum-Connes conjecture. 
\begin{theorem} \label{bcqut}
Let $ A $ be a finite dimensional $ C^* $-algebra of dimension $ \dim(A) \geq 4 $, and let $ G $ be the quantum automorphism group 
of $ A $ with respect to a $ \delta $-form on $ A $. Then we have $ KK^G = \bra \T_G \ket $.  
In particular, the dual of $ G $ is $ K $-amenable. 
\end{theorem} 
\proof For the first claim it is enough to check the corresponding assertion for $ H = SO_q(3) $. 
In this case it follows from the strong Baum-Connes property of $ K = SU_q(2) $, see \cite{Voigtbcfo}. More precisely, 
let $ B \in KK^{\check{H}} $ and consider the induction 
functor $ KK^{\check{H}} \rightarrow KK^{\check{K}} $. According to \cite{Voigtbcfo}, the algebra $ \ind_{\check{H}}^{\check{K}}(B) $ is contained 
in the localising subcategory of $ KK^{\check{K}} $ generated by all algebras of the form $ C_0(\check{K}) \otimes C $, where 
$ C $ is any $ C^* $-algebra, with the coaction given by comultiplication on the first tensor factor. Moreover, due to lemma \ref{soq3morita}, 
the $ \check{H} $-$ C^* $-algebra $ C_0(\check{K}) $ decomposes as 
$$ 
C_0(\check{K}) \cong (SO_q(3) \ltimes \mathbb{C}) \oplus (SO_q(3) \ltimes \KH(V(1/2)))
$$ 
in $ KK^{\check{H}} $. Therefore it suffices to observe that $ B $ is a direct summand 
in $ \res^{\check{K}}_{\check{H}} \ind_{\check{H}}^{\check{K}}(B) $ as an $ \check{H} $-algebra. Using the concrete 
description of induced $ C^* $-algebras given in \cite{VVfreeu}, this in turn follows by considering 
the coaction of $ B $, viewed as a $ * $-homomorphism $ \beta: B \rightarrow \ind_{\check{H}}^{\check{K}}(B) \subset M(C_0(\hat{K}) \otimes B) $, 
and the map $ \check{\epsilon} \otimes \id: \ind_{\check{H}}^{\check{K}}(B) \rightarrow B $ where $ \check{\epsilon} $ denotes 
the counit of $ C_0(\check{K}) $. \\ 
The $ K $-amenability of $ \check{G} $ follows now immediately using lemma \ref{properamenable}, compare the analogous argument 
in \cite{Voigtbcfo}. \qed \\
If $ C $ is an $ SO_q(3) $-$ C^* $-algebra we write $ \KH(V(1/2)) \otimes C $ 
for the $ SO_q(3) $-algebra obtained by considering the $ SU_q(2) $-$ C^* $-algebra $ \KH(V(1/2)) \otimes C $ and observing 
that the coaction on $ \KH(V(1/2)) \otimes C $ takes values in $ C(SO_q(3)) \otimes \KH(V(1/2)) \otimes C $. 
\begin{lemma} \label{Khalfshifting}
For any $ SO_q(3) $-$ C^* $-algebra $ B $ we have an $ SO_q(3) $-equivariant Morita equivalence 
$$ 
\KH(V(1/2)) \otimes \KH(V(1/2)) \otimes B \sim_M B. 
$$ 
Moreover, there exists a natural isomorphism 
$$
KK^{SO_q(3)}(\KH(V(1/2))  \otimes B, \KH(V(1/2)) \otimes C) \cong KK^{SO_q(3)}(B, C)
$$
for all $ SO_q(3) $-$ C^* $-algebras $ B, C $. 
\end{lemma} 
\proof For the first claim observe that 
$$ 
\KH(V(1/2)) \otimes \KH(V(1/2)) \otimes B \cong \KH(V(1/2) \otimes V(1/2)) \otimes B 
$$ 
is $ SO_q(3) $-equivariantly Morita equivalent to $ B $ since $ V(1/2) \otimes V(1/2) \cong V(0) \oplus V(1) $ 
is an honest representation of $ SO_q(3) $. \\ 
The second part of the lemma follows from this. More precisely, the functor $ F: SO_q(3) \Alg \rightarrow KK^{SO_q(3)} $ given 
by $ F(B) = \KH(V(1/2)) \otimes B $ is homotopy invariant, stable, and split exact. By the universal property of 
equivariant $ KK $-theory \cite{NVpoincare}, it therefore induces a functor $ f: KK^{SO_q(3)} \rightarrow KK^{SO_q(3)} $. Alternatively, 
one can also construct $ f $ directly on the level of Kasparov cycles. 
Using the above Morita equivalence, one checks that $ f^2 $ is naturally isomorphic to 
the identity. In particular, $ f $ is a natural isomorphism. \qed \\ 
We need some basic calculations. Let us write 
$$ 
A = M_{k_1}(\mathbb{C}) \oplus \cdots \oplus M_{k_n}(\mathbb{C})
$$ 
and view it as a $ G $-$ C^* $-algebra with the defining action of $ G = \Qut(A, \omega) $. 
Moreover we write write $ \mathbb{C} $ for the trivial $ G $-$ C^* $-algebra. 
As indicated above, we may identify $ \Irr(G) \cong \mathbb{N}_0 $ and $ R(G) = \mathbb{Z}[t] $, 
where $ t $ corresponds to the underlying representation of $ A $, that is, $ t = V(0) + V(1) $ under the identification 
with $ R(SO_q(3)) = K_*(C^*(SO_q(3))) $. Similarly, we let $ R^\omega(G) = K_*(C^*_\omega(SO_q(3))) $, 
and identify $ R^\omega(G) = t^{1/2} \mathbb{Z}[t] $ as an $ R(G) $-module. Here the generator $ t^{1/2} $ corresponds to $ V(1/2) $. \\ 
The representation ring $ R(G) $ acts on the groups $ KK^G(B,C) $ in a natural way, compare \cite{Meyerhomalg2}. We have 
to identify these $ R(G) $-modules in a few special cases. 
\begin{lemma} \label{repringcomputations}
Let $ G = \Qut(A, \omega) $ as above. Then we have isomorphisms 
\begin{align*}
KK^G(\mathbb{C}, C^\red(G)) &\cong KK(\mathbb{C}, \mathbb{C}) = \mathbb{Z} \\ 
KK^G(A, C^\red(G)) &\cong KK(A, \mathbb{C}) = \mathbb{Z}^n 
\end{align*}
and 
\begin{align*}
KK^G(\mathbb{C}, \mathbb{C}) \cong R(&G) \cong KK^G(A, A) \\ 
KK^G(A, \mathbb{C}) \cong R^\omega(&G) \cong KK^G(\mathbb{C}, A)
\end{align*}
of $ R(G) $-modules. 
\end{lemma}
\proof In the same way as in the proof of proposition 4.7 in \cite{NVpoincare} one checks the first two isomorphisms. 
Indeed, although the quantum group $ G $ will typically fail to be coamenable, the arguments given there essentially carry over because 
$ \check{G} $ is $ K $-amenable, see theorem \ref{bcqut}. The remaining isomorphisms follow from lemma \ref{Khalfshifting}. \\ 
The $ R(G) $-module structures on the various Kasparov groups involving $ A $ and $ \mathbb{C} $ correspond to the obvious ones 
on $ R(G) $ and $ R^\omega(G) $. \qed \\ 
We are now ready to prove our main theorem. 
\begin{theorem} \label{Kqaut}
Let 
$$ 
A = M_{k_1}(\mathbb{C}) \oplus \cdots \oplus M_{k_n}(\mathbb{C})
$$ 
be a finite dimensional $ C^* $-algebra of dimension at least $ 4 $, equipped with a $ \delta $-form $ \omega $ for some $ \delta > 0 $. 
Then 
$$
K_0(C(G)) = \mathbb{Z}^{(n - 1)^2 + 1} \oplus \mathbb{Z}_d^{2n - 1}, \qquad K_1(C(G)) = \mathbb{Z}, 
$$
where $ \mathbb{Z}_d $ is the finite cyclic group of order $ d = \gcd(k_1, \dots, k_n) $,
the greatest common divisor of the numbers $ k_1, \dots, k_n $. 
\end{theorem}
\proof Notice that according to theorem \ref{bcqut} the dual of the quantum automorphism group $ G = \Qut(A, \omega) $ is $ K $-amenable, 
so that we do not have to distinguish between $ C^\max(G) $ and $ C^\red(G) $ as far as $ K $-theory is concerned. \\ 
In order to compute the $ K $-groups $ K_*(C^\max(G)) \cong K_*(C^\red(G)) $ we have to write down a suitable resolution of the trivial 
action on $ \mathbb{C} $ in $ KK^{\check{G}} $. It is slightly easier to work in the compact picture, so that we shall construct a 
resolution of $ C^\red(G) $ in $ KK^G $. \\ 
The corresponding homological algebra can be expressed using the framework of homological ideals, see \cite{MNhomalg1}. 
In order to do this it is convenient to first pass from $ KK^G $ to $ KK^H $ for $ H = SO_q(3) $ monoidally equivalent to $ G $. 
We define an ideal $ \mathfrak{J}_H \subset KK^H $ by taking the intersection $ \mathfrak{J}_H = \ker(F_0) \cap \ker(F_1) $ 
where $ F_j: KK^H \rightarrow KK $ are the functors given by 
$$
F_0(C) = H \ltimes C, \qquad F_1(C) = H \ltimes (\KH(V(1/2)) \otimes C), 
$$
respectively. It is straightforward to check that $ \mathfrak{J}_H $ is a stable homological ideal in $ KK^H $, 
and we let $ \mathfrak{J}_G $ be the corresponding ideal in $ KK^G $. The general 
machinery in \cite{MNhomalg1}, \cite{Meyerhomalg2} now allows us to study $ \mathfrak{J}_G $-projective resolutions in $ KK^G $. \\ 
Recall from lemma \ref{quttorsion} that the defining coaction on $ A $ and the trivial coaction on $ \mathbb{C} $ 
are the only torsion coactions of $ \check{G} $ up to equivariant Morita equivalence. It is straightforward to check that 
both these algebras are $ \mathfrak{J}_G $-projective when viewed as objects in $ KK^G $. \\ 
Let us consider the diagram $ C_\bullet $ given by 
$$
\xymatrix{
0 \ar@{->}[r] & C_1 \ar@{->}[r]^{d_1} & C_0 \ar@{->}[r]^{d_0} & C^\red(G) \ar@{->}[r] & 0
}
$$
in $ KK^G $ where 
\begin{align*}
C_1 &= \mathbb{C}^{\oplus n} \oplus A \\ 
C_0 &= A^{\oplus n} \oplus \mathbb{C}
\end{align*}
and the arrows are defined as follows. The morphism $ d_0 $ is 
$$
d_0 = \epsilon_1 \oplus \cdots \oplus \epsilon_n \oplus u
$$
where $ \epsilon_j $ corresponds to the canonical basis vector $ e_j $ in $ \mathbb{Z}^n = KK^G(A, C^\red(G)) $ and 
$ u: \mathbb{C} \rightarrow C^\red(G) $ is the unit homomorphism. The morphism $ d_1 $ is 
$$
d_1 = \begin{pmatrix} 
t^{1/2} & 0 & \hdots & & -k_1 \\
0 & \ddots & & & \vdots \\ 
\vdots & & \ddots & & \\ 
 & & & t^{1/2} & -k_n \\ 
-k_1 & \cdots & & -k_n & t^{1/2}, 
\end{pmatrix} 
$$
where we use the identifications obtained in lemma \ref{repringcomputations}. \\ 
In order to show that $ C_\bullet $ is $ \mathfrak{J} $-exact, it suffices to check that the sequences $ KK^G(\mathbb{C}, C_\bullet) $
and $ KK^G(A, C_\bullet) $ are both exact. \\ 
Let us compute the induced maps in the sequence $ KK^G(\mathbb{C}, C_\bullet) $. For $ d_0 $ we obtain 
$$
d^\mathbb{C}_0 = 
\begin{pmatrix} 
k_1 &
\cdots &
k_n & 
1 
\end{pmatrix}, 
$$
viewed as the $ R(G) $-linear map acting on the free $ R(G) $-module $ KK^G(\mathbb{C}, C_0) \cong R^\omega(G)^{\oplus n} \oplus R(G) \cong R(G)^{\oplus n + 1}$. 
Indeed, the generator $ t^{1/2} $ in the $ j $-th copy of this module corresponds to the unit map $ \mathbb{C} \rightarrow A $, 
and composing with $ \epsilon_j $ yields $ k_j \in \mathbb{Z} = KK^G(\mathbb{C}, C^\red(G)$. 
Moreover, the generator $ 1 $ in the last copy of $ R(G) $ is clearly mapped to $ 1 \in \mathbb{Z} $. \\ 
The morphism $ d_1^\mathbb{C} $ becomes  
$$
d^\mathbb{C}_1 = \begin{pmatrix} 
1 & 0 & \hdots & & -k_1 \\
0 & \ddots & & & \vdots \\ 
\vdots & & \ddots & & \\ 
  & & & 1 & -k_n\\ 
-k_1 & \cdots & & -k_n & t 
\end{pmatrix}. 
$$
It is straightforward to check that $ KK^G(\mathbb{C}, C_\bullet) $ is an exact complex. \\ 
Similarly, for the induced maps in $ KK^G(A, C_\bullet) $ we obtain 
$$
d^A_0 = 
\begin{pmatrix} 
1 & 0 & \hdots & & & k_1 \\
0 & \ddots & & & & \vdots \\ 
\vdots & & \ddots & & & \vdots \\ 
  & & & \ddots & & \vdots \\ 
0 & \cdots & & 0 & 1 & k_n 
\end{pmatrix} 
$$
and 
$$
d^A_1 = \begin{pmatrix} 
t & 0 & \hdots & & -k_1 \\
0 & \ddots & & & \vdots \\ 
\vdots & & \ddots & & \vdots \\ 
  & & & t & -k_n \\ 
-k_1 & \cdots & & -k_n & 1
\end{pmatrix}.  
$$
The resulting diagram is an exact complex as well, but this is slightly more subtle. Indeed, one has to be careful to take into account 
the correct $ R(G) $-module structures when identifying the action of $ t $. 
With this in mind, in order to check $ \ker(d^A_0) = \im(d^A_1) $ 
one has to inductively reduce every element in $ \ker(d^A_0) $ to an $ n $-tuple of constant polynomials modulo $ \im(d^A_1) $, and such 
elements are in the image of $ d_0^A $. The surjectivity of $ d^A_0 $ and the injectivity of $ d^A_1 $ are easy. \\ 
It follows that $ C_\bullet $ is indeed a $ \mathfrak{J} $-projective resolution, and to compute the 
$ K $-theory of $ C(G) $ it suffices to compute kernel and cokernel of the map 
$ \partial: \mathbb{Z}^n \oplus \mathbb{Z}^n = K_0(\mathbb{C}^{\oplus n} \oplus A) \rightarrow 
K_0(A^{\oplus n} \oplus \mathbb{C}) = (\mathbb{Z}^n)^n \oplus \mathbb{Z} $ 
given by 
$$
\partial = \begin{pmatrix} 
{\bf k}^T & 0 & \hdots & & -k_1 {\bf 1} \\
0 & \ddots & & & \vdots \\ 
\vdots & & \ddots & & \vdots \\ 
  & & & {\bf k}^T & -k_n {\bf 1} \\ 
-k_1 & \cdots & & -k_n & {\bf k}
\end{pmatrix}.  
$$
Here $ {\bf k}^T $ is the transpose of $ (k_1, \dots, k_n) = {\bf k} $ and $ {\bf 1} $ is the identity matrix in $ M_n(\mathbb{Z}) $. \\
Let us compute $ \ker(\partial) $. Inspecting the $ rn + 1 $th rows of $ \partial $ for $ r = 0, \dots, n - 1 $ we see that an 
element of $ \ker(\partial) $ is necessarily of the form $ (a_1, \dots, a_n, a_1, \dots, a_n) $. Moreover, the first $ n $ rows give the relations 
$ k_i a_1 = k_1 a_i $, the next $ n $ rows give $ k_i a_2 = k_2 a_i $, and so on. That is, we obtain 
$$
k_i a_j = k_j a_i 
$$ 
for all $ 1 \leq i,j \leq n $. Notice in particular that $ a_2, \dots, a_n $ are uniquely determined by $ a_1 \in \mathbb{Z} $. \\
The general solution to these equations is $ (a_1, \dots, a_n) = (m k_1/d, \dots, m k_n/d) $ where $ m \in \mathbb{Z} $ and 
$ d = gcd(k_1, \dots, k_n) $ is the greatest common divisor of $ k_1, \dots, k_n $. In particular, we conclude 
$$ 
\ker(\partial) = \mathbb{Z}. 
$$ 
Let us now compute $ \coker(\partial) $. Using elementary row operations we can transform $ \partial $ to 
$$
\partial_1 = \begin{pmatrix} 
{\bf k}^T & 0 & \hdots & & -k_1 {\bf 1} \\
0 & \ddots & & & \vdots \\ 
\vdots & & \ddots & & \vdots \\ 
  & & & {\bf k}^T & -k_n {\bf 1} \\ 
0 & \cdots & & 0 & {\bf 0}
\end{pmatrix}
$$
We thus obtain a direct summand $ \mathbb{Z} $ in $ \coker(\partial) $, and we may restrict to the matrix obtained 
by deleting the last row from $ \partial_1 $. Performing elementary row and column operations we may reduce the 
resulting matrix to 
$$
\partial_2 = \begin{pmatrix} 
{\bf v}^T & 0 & \hdots & & -k_1 {\bf 1} \\
0 & \ddots & & & \vdots \\ 
\vdots & & \ddots & & \vdots \\ 
  & & & {\bf v}^T & -k_n {\bf 1}
\end{pmatrix}
$$
where $ {\bf v} = (d, 0, \dots, 0) $, and we recall that $ d = gcd(k_1, \dots, k_n) $. 
Simplifying the right hand side of this matrix further leads to a diagonal matrix with $ n + (n - 1) $ entries $ d $, 
and all remaing entries zero. Hence the final result is 
$$ 
\coker(\partial) = \mathbb{Z}^{(n - 1)^2 + 1} \oplus \mathbb{Z}_d^{2n - 1}
$$ 
as claimed. \qed \\ 
Let us remark that the case $ n = 1 $ of theorem \ref{Kqaut} was already discussed in \cite{Voigtbcsuq2}. 
At the opposite extreme $ k_1 = \cdots = k_n = 1 $, theorem \ref{Kqaut} implies the following result. 
\begin{cor} \label{snplus}
Let $ n \geq 4 $. Then the quantum permutation group $ S_n^+ $ is $ K $-amenable, and the $ K $-theory is given by 
$$
K_0(C(S_n^+)) = \mathbb{Z}^{n^2 - 2n + 2}, \qquad K_1(C(S_n^+)) = \mathbb{Z}.  
$$
Generators in degree zero are given by the projections $ 1 $ and $ u_{ij} $ for $ 1 \leq i,j \leq n - 1 $. 
\end{cor} 
\proof It remains only to verify the claim regarding generators of $ K_0(C(S_n^+)) $. 
For this it suffices to consider the images of the generating projections $ u_{ij} \in C^\max(S_n^+) $ in $ K_0(C(S_n)) $, and notice that 
they span a copy of $ \mathbb{Z}^{n^2 - 2n + 2} $ inside $ K_0(C(S_n)) = \mathbb{Z}^{n!} $. Essentially, 
in each row and column of the matrix $ u = (u_{ij}) $ the last entry is determined by the remaining $ n - 1 $ entries, 
with no further relations. This accounts for $ (n - 1)^2 = n^2 - 2n + 1 $ generators, and 
the missing generator is the class of the unit. \qed \\
By mapping $ C(S_n^+) $ to $ C(S_4^+) $ and using theorem 5.2 in \cite{BBfourpoints}, it is not hard to check that the defining 
unitary $ u \in M_n(C(S_n^+)) $ yields a nonzero class $ [u] \in K_1(C(S_n^+)) = \mathbb{Z} $, and that $ [u] $ is of the form 
$ k x $ where $ x $ is a generator and $ k \in \mathbb{N} $ is at most $ 8 $. 
However, to actually identify the generator of $ K_1(C(S_n^+)) $ would require more work. \\ 
Let us point our that corollary \ref{snplus} shows in particular that the $ K $-theory of $ C(S_4^+) = C(SO_{-1}(3)) $ differs significantly 
from the $ K $-theory of $ SO(3) $. \\ 
Using the explicit structure of $ C^\red(S_n^+) $ for $ n = 1,2,3 $ and corollary \ref{snplus}, we can distinguish 
the reduced $ C^* $-algebras $ C^\red(S_n^+) $ for different values of $ n $. 
\begin{cor} 
Let $ m, n \in \mathbb{N} $. Then $ C^\red(S_m^+) \cong C^\red(S_n^+) $ iff $ m = n $. 
\end{cor} 
Of course, this result holds for the maximal $ C^* $-algebras as well. Notice that the maximal $ C^* $-algebras can 
already be distinguished by comparing their abelianisations, a method which does not work for the reduced $ C^* $-algebras.

\bibliographystyle{plain}

\bibliography{cvoigt}

\end{document}